\documentclass{gtart_h}

\def\ifplaintex{\expandafter\ifx\csname documentclass\endcsname\relax}

\def\gtp{{\mathsurround=0pt\it $\cal G\mskip-2mu$eometry \&\ 
$\cal T\!\!$opology $\cal P\!$ublications}}  

\def\recd{{\small Received:\qua\receiveddate\ifx\reviseddate\relax
\else\qquad Revised:\qua\reviseddate\fi\par}} 

\def\lognumber#1{\def\thelognumber{#1}}
\def\volumenumber#1{\def\thevolumenumber{#1}}
\def\volumeyear#1{\def\thevolumeyear{#1}}
\def\papernumber#1{\def\thepapernumber{#1}}
\def\pagenumbers#1#2{\def\startpage{#1}\def\finishpage{#2}}
\def\published#1{\def\publishdate{#1}}

\def\received#1{\def\receiveddate{#1}}

\def\accepted#1{\def\accepteddate{#1}}

\long\def\asciiabstract#1{\long\def\theasciiabstract{#1}}
\def\asciikeywords#1{\def\theasciikeywords{#1}}

\let\\\par\let\thelognumber\relax\let\thevolumenumber\relax
\let\thepapernumber\relax\let\thevolumeyear\relax\let\startpage\relax
\let\finishpage\relax\let\publishdate\relax\let\receiveddate\relax
\let\reviseddate\relax\let\accepteddate\relax\let\theasciititle\relax
\let\theasciiauthors\relax
\let\theasciiabstract\relax\let\theasciikeywords\relax

\let\theasciiemail\relax

\ifplaintex
\font\logobig=cmssbx10 scaled 3836
\font\logomed=cmssbx10 scaled 2557
\else
\font\logobig=cmssbx10 scaled 4200
\font\logomed=cmssbx10 scaled 2800
\fi

\long\def\makeagttitle{   
\count0=\startpage
\agt\hfill      
\hbox to 45truept{\vbox to 0pt{\vglue -13truept{\logomed A\kern -.37em{\logobig 
T}\kern -.38em G}\vss}\hss}
\break
{\small Volume \thevolumenumber\ (\thevolumeyear)
\startpage--\finishpage\nl
Published: \publishdate}

\vglue .25truein

{\parskip=0pt\leftskip 0pt plus
1fil\def\\{\par\smallskip}{\Large\bf\thetitle}\par\medskip} \vglue
0.05truein

{\parskip=0pt\leftskip 0pt plus 1fil\def\\{\par}{\sc\theauthors}
\par\medskip}%
 
\vglue 0.03truein 

{\small\leftskip 25truept\rightskip 25truept{\bf Abstract}\stdspace\theabstract

{\bf AMS Classification}\stdspace\theprimaryclass
\ifx\thesecondaryclass\relax\else; \thesecondaryclass\fi\par
{\bf Keywords}\stdspace \thekeywords\par}\vglue 7truept

}   

\ifplaintex
\font\phead=cmsl9 scaled 950
\font\pnum=cmbx10 scaled 913
\font\pfoot=cmsl9 scaled 950
\headline{\vbox to 0pt{\vskip -4.5mm\line{\small\phead\ifnum
\count0=\startpage ISSN 1472-2739 (on-line) 1472-2747 (printed)
\hfill {\pnum\folio}\else\ifodd\count0\def\\{ }%
\ifx\theshorttitle\relax\thetitle\else\theshorttitle\fi\hfill{\pnum\folio}
\else\def\\{ and }{\pnum\folio}\hfill\ifx\theshortauthors\relax\theauthors
\else\theshortauthors\fi\fi\fi}\vss}}
\footline{\vbox to 0pt{\vglue 0mm\line{\small\pfoot\ifnum\count0=\startpage
\copyright\ \gtp\hfill\else
\agt, Volume \thevolumenumber\ (\thevolumeyear)\hfill\fi}\vss}}
\else
\font=cmsl9 scaled 1050
\font\lnum=cmbx10 
\font=cmsl9 scaled 1050
\makeatletter
\def\@oddhead{{\small\lhead\ifnum\count0=\startpage ISSN 1472-2739 
(on-line) 1472-2747 (printed)\hfill {\lnum\number\count0}\else\ifodd\count0
\def\\{ }\ifx\theshorttitle\relax \thetitle \else\theshorttitle\fi\hfill
{\lnum\number\count0}\else\def\\{ and }{\lnum\number\count0}
\hfill\ifx\theshortauthors\relax 
\theauthors\else\theshortauthors\fi\fi\fi}}\def\@evenhead{\@oddhead}
\def\@oddfoot{\small\lfoot\ifnum\count0=\startpage\copyright\ \gtp\hfill\else
\agt, Volume \thevolumenumber\ (\thevolumeyear)\hfill\fi}
\def\@evenfoot{\@oddfoot}
\makeatother
\fi
\let\maketitlepage\makeagttitle
\let\makeshorttitle\maketitlepage
\let\maketitle\maketitlepage

\newwrite\gtoutfile
\long\gdef\makeheadfile{  
{\def\\{, }\def\s{ }
\immediate\openout\gtoutfile head.xxx
\immediate\write\gtoutfile{Proxy-for: \ifx\theasciiauthors\relax
\theauthors\else\theasciiauthors\fi\s<\ifx\theasciiemail\relax\theemail\else\theasciiemail\fi>}
\immediate\write\gtoutfile{\noexpand\\}
\immediate\write\gtoutfile{Authors: \ifx\theasciiauthors\relax
\theauthors\else\theasciiauthors\fi}
{\def\\{ }\immediate\write\gtoutfile{Title: \ifx\theasciititle\relax
\thetitle\else\theasciititle\fi}}
\immediate\write\gtoutfile{Subj-class: GT or SG, GR etc}
\immediate\write\gtoutfile{MSC-class: \theprimaryclass\ifx\thesecondaryclass\relax\else, \thesecondaryclass\fi}
\immediate\write\gtoutfile{Journal-ref: Algebraic and Geometric Topology \thevolumenumber\s
(\thevolumeyear) \startpage-\finishpage}
\immediate\write\gtoutfile{Comments: Published by Algebraic and
Geometric Topology at}
\immediate\write\gtoutfile{\s\s\s  http://www.maths.warwick.ac.uk/agt/AGTVol\thevolumenumber/agt-\thevolumenumber-\thepapernumber.abs.html}
\immediate\write\gtoutfile{\noexpand\\}
\immediate\write\gtoutfile{}
\ifx\theasciiabstract\relax
\immediate\write\gtoutfile{\theabstract}\else
\immediate\write\gtoutfile{\theasciiabstract}\fi
\immediate\write\gtoutfile{}
\immediate\write\gtoutfile{\noexpand\\}
\immediate\write\gtoutfile{}
\immediate\closeout\gtoutfile}}  

\def\maketitlepage{\makeagttitle\makeheadfile}
\let\makeshorttitle\maketitlepage
\let\maketitle\maketitlepage

\lognumber{16}
\volumenumber{4}
\volumeyear{2004}
\papernumber{16}
\published{6 May 2004}
\pagenumbers{297}{309}
\received{31 August 2003}
\accepted{21 March 2004}

\usepackage{amsmath,amssymb,amscd,graphicx}

\newcommand{\R}{\mathbb{R}}
\newcommand{\cat}{\textsc{CAT(0)}}

\newtheorem{thm}{Theorem}[section]
\newtheorem*{thmn}{Theorem}
\newtheorem{lem}[thm]{Lemma}
\newtheorem{cor}[thm]{Corollary}

\newtheorem{prop}[thm]{Proposition}

\theoremstyle{definition}
\newtheorem{defn}[thm]{Definition}

\newtheorem{ex}[thm]{Example}

\begin{document}

\title{Cubulating spaces with walls}
\author{Bogdan Nica}
\address{Department of Mathematics, McGill 
University\\Montreal, Canada H3A 2K6}
\email{bnica@math.mcgill.ca}

\begin{abstract}
We describe a correspondence between spaces with walls and $\cat$
cube complexes.
\end{abstract}
\asciiabstract{%
We describe a correspondence between spaces with walls and CAT(0)
cube complexes.}


\asciikeywords{Space with walls, Median graph, CAT(0)
cube complex}

\maketitle

\section{Introduction}
\noindent The elegant notion of a space with walls was introduced
by Haglund and Paulin \cite{HP98}. Prototypical examples of spaces
with walls are $\cat$ cube complexes, introduced by Gromov in
\cite{Gro87}. The purpose of this note is to observe that every
space with walls has a canonical embedding in a $\cat$ cube
complex and, consequently, a group action on a space with walls
extends naturally to a group action on a $\cat$ cube complex. The
usefulness of this result is that spaces with
walls are often easily identifiable by geometric reasons.

The cubulation of a space with walls, as we call it, is an
abstract version of a construction introduced by Sageev
\cite{Sag95} for the purpose of relating multi-ended pairs of
groups to essential actions on $\cat$ cube complexes. Sageev's
construction is further explored by Niblo and Roller in
\cite{NRo98}, where an essential group action on a $\cat$ cube
complex is shown to imply the failure of Kazhdan's property (T)
(see also \cite{NRe97}). Roller's detailed study \cite{Ro}
formulates Sageev's construction in the language of median
algebras (see also \cite{Ger98}). Finally, a version of Sageev's
construction, where a $\cat$ cube complex arises from a system of
halfspaces in a complex, is considered by Niblo and Reeves
\cite{NRe03}, for Coxeter groups, and by Wise \cite{Wi}, for
certain small cancellation groups.

Some of the papers cited above (\cite{Sag95}, \cite{NRo98},
\cite{NRe03}) take the point of view that the primitive data for
constructing a $\cat$ cube complex is a partially ordered set with
an order-reversing involution, with certain discreteness and
nesting assumptions \`{a} la Dunwoody, which is to become the
system of halfspaces in the cube complex. However, a space with
walls comes in handy when a suitable connected component needs to
be specified.

The cubulation of a space with walls comprises two steps,
according to the following scheme:
\begin{displaymath}
\begin{CD}
\textrm{ space with walls } @>>> \textrm{ median graph } @>>>
\cat\textrm{ cube complex}
\end{CD}
\end{displaymath}
In the first step, which is our main objective, a space with walls
$X$ is embedded in a median graph $\mathcal{C}^1(X)$, called the
``1--cubulation of the space with walls $X$''. The second step is
based on the fact that any median graph is the 1--skeleton of a
unique $\cat$ cube complex. Explicitly, the step from the median
graph $\mathcal{C}^1(X)$ to a $\cat$ cube complex
$\mathcal{C}(X)$ consists of ``filling in'' isometric copies of
euclidean cubes by inductively adding an $n$--dimensional cube
whenever its $(n-1)$--skeleton is present; see \cite[\S3]{Sag95},
\cite[\S6]{Che00}, \cite[Thm.10.3]{Ro}, \cite[\S5]{Wi}.

Our result is:
\begin{thmn}
Let $X$ be a space with walls. There exists an injective morphism
of spaces with walls
\begin{displaymath}
\begin{CD}
X @>\sigma>> \mathcal{C}^1(X)
\end{CD}
\end{displaymath}
where $\mathcal{C}^1(X)$ is a connected median graph, and
$\sigma(X)$ ``spans'' $\mathcal{C}^1(X)$, in the sense that no
proper subgraph of $\mathcal{C}^1(X)$ containing $\sigma(X)$ is
median.
\end{thmn}

Spaces with walls and their morphisms are defined in Section 2.
Median graphs are defined in Section 3, where we prove that they
are spaces with walls. The main construction is presented in
Section 4, where we also show that any group action on a space
with walls has a unique extension to a group action on its
1--cubulation.

Finally, we would like to draw the reader's attention to a
different account of the cubulation procedure, described
independently by Chatterji and Niblo \cite{CN}.

\rk{Acknowledgements} I am grateful to Dani Wise for suggesting
the problem, as well as for constructive comments that have
improved the content of this paper. This work has been supported
by NSERC.

\section{Spaces with walls}
\noindent We recall the definition of a space with walls. We also
introduce the natural notion of morphism of spaces with walls.

\begin{defn}
Let $X$ be a set. A \emph{wall} in $X$ is a partition of $X$ into
2 subsets called \emph{halfspaces}. We say that $X$ is a
\emph{space with walls} if $X$ is endowed with a collection of
walls, containing the trivial wall $\{\emptyset, X\}$, and so
that any two distinct points are separated by a finite, non-zero
number of walls. Note that a wall separates two distinct points
$x,y\in X$ if $x$ belongs to one of the halfspaces determined by
the wall, while $y$ belongs to the other halfspace.

A \emph{morphism} of spaces with walls is a map
$f\co\negthinspace X\negthinspace\to \negthinspace X'$ between
spaces with walls with the property that $f^{-1}(A')$ is a
halfspace of $X$ for each halfspace $A'$ of $X'$.
\end{defn}

A minor difference between the original definition from
\cite{HP98} and the definition given above is that we insist on
the presence of the trivial wall. This modification is needed for
a morphism of spaces with walls to be well-defined. Another
reason is that halfspaces arise naturally in the presence of an
underlying convexity structure: a halfspace is a convex set whose
complement is convex. In such a context, the trivial wall is
always present.

A space with walls $X$ becomes a metric space by defining the
distance between two points to be the number of walls separating
them: $d_w(x,y)=\big|W(x,y)\big|$, where $W(x,y)$ denotes the set
of walls separating $x$ and $y$. For the triangle inequality
observe that, given $x,y,z$, a wall separating $x$ and $y$ has to
separate either $x$ from $z$, or $z$ from $y$, i.e.,
$W(x,y)\subseteq W(x,z)\cup W(z,y)$. We obtain:
$d_w(x,y)=\big|W(x,y)\big|\leq \big|W(x,z)\big|+
\big|W(z,y)\big|=d_w(x,z)+d_w(z,y)$. We call $d_w$ the \emph{wall
metric} on $X$.

A group acts on a space with walls $X$ by permuting the walls.
Consequently, it acts by isometries on $(X,d_w)$.

\section{Median graphs}

\noindent Median graphs are well documented in graph-theoretic
literature. See \cite[\S4]{Che00} for a list of papers on median
graphs. Unless otherwise specified, graphs are henceforth assumed
to be connected, simplicial, in the sense that they have no loops
or multiple edges, and equipped with the path metric. The
\emph{geodesic interval} $[x,y]$ determined by the vertices $x$
and $y$ is the collection of vertices lying on a shortest path
from $x$ to $y$.

\begin{defn}
A graph is \emph{median} if, for each triple of vertices $x$, $y$,
$z$, the geodesic intervals $[x,y]$, $[y,z]$, $[z,x]$ have a
unique common point.
\end{defn}

\noindent Trees are elementary examples of median graphs. The
1--skeleton of the square tiling of the plane is median whereas
the 1--skeletons of the hexagonal and triangular tilings are not.

\begin{lem}
In a median graph the geodesic intervals are finite.
\end{lem}
\begin{proof} Induction on the distance between vertices. Consider two vertices $x,y$
and let $Z$ denote the collection of those neighbors of $y$ that
are in $[x,y]$. As $[x,y]\setminus \{y\}\subseteq \cup_{z\in Z}
[x,z]$ with each $[x,z]$ finite by the induction hypothesis, it
suffices to show that $Z$ is finite.

Nodes in $Z$ are pairwise distance 2 apart. Fix $a\in Z$ and let
$m(z)=m(x,a,z)$ for every $z\in Z$. See Figure~\ref{fig:Kite}.

\begin{figure}[ht!]\centering
\includegraphics[width=.5\textwidth]{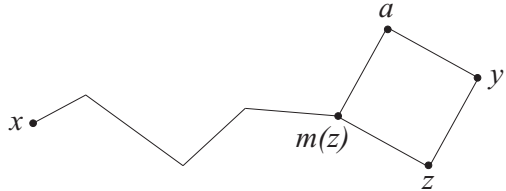}
\caption{Finiteness of geodesic intervals in a median graph
\label{fig:Kite}}
\end{figure}

The mapping $z\mapsto m(z)$ from $Z$ to $[x,a]$ is injective: if
$m(z)=m(z')=m$ for distinct $z,z'\in Z$ then both $m$ and $y$ are
medians for the triple $a,z,z'$. The finiteness of $Z$ follows
now from the finiteness of $[x,a]$.
\end{proof}

\noindent Median graphs are instances of median algebras, which
can be described as interval structures that enjoy a tripod-like
condition. Median algebras can also be defined as sets with a
ternary operation $x,y,z\mapsto m(x,y,z)$, called the median
operation, that satisfies certain axioms. The latter point of
view is adopted in \cite{BH83} and \cite{Ro}, which we suggest as
supplements to our brief presentation.

\begin{defn}
A \emph{median algebra} is a set $X$ with an interval assignment
$(x,y)\mapsto [x,y]$, mapping pairs of points in $X$ to subsets of
$X$, so that for all $x,y,z\in X$ the following are satisfied:

$\bullet$ $[x,x]=\{x\}$

$\bullet$ if $z\in[x,y]$ then $[x,z]\subseteq [x,y]$

$\bullet$ $[x,y]$, $[y,z]$, $[z,x]$ have a unique common point,
denoted $m(x,y,z)$

\noindent A subset $A\subseteq X$ is \emph{convex} if
$[x,y]\subseteq A$ for all $x,y\in A$. A subset $A\subseteq X$ is
a \emph{halfspace} if both $A$ and $A^c=X\setminus A$ are convex.

A \emph{morphism} of median algebras is a map $f\co \negthinspace
X\negthinspace\to \negthinspace X'$ between median algebras that
is ``betweenness preserving'', in the sense that
$f\big([x,y]\big)\subseteq \big[f(x),f(y)\big]$ for all $x,y\in
X$.
\end{defn}

\begin{defn}
A median algebra is \emph{discrete} if every interval is finite.
\end{defn}

\begin{ex}\label{boolean}(Boolean median algebra)\qua Any power set $\mathcal{P}(X)$ is a
median algebra under the interval assignment
\begin{displaymath}
(A,B)\mapsto [A,B]=\{C: A\cap B\subseteq C\subseteq A\cup B\}
\end{displaymath}
with $(A\cap B)\cup(B\cap C)\cup(C\cap A)=(A\cup B)\cap(B\cup
C)\cap(C\cup A)$ being the boolean median of $A$, $B$, $C$. The
significance of this example is that any median algebra is
isomorphic to a subalgebra of a boolean median algebra.

In general, a boolean median algebra $\mathcal{P}(X)$ is far from
discrete. We say that $A,B\subseteq X$ are \emph{almost equal} if
their symmetric difference $A\triangle B=(A\setminus B)
\cup(B\setminus A)$ is finite; this defines an equivalence
relation on $\mathcal{P}(X)$. An almost equality class of
$\mathcal{P}(X)$ is a discrete median algebra.
\end{ex}

A crucial feature of halfspaces in median algebras is that they
separate disjoint convex sets, in the sense that for all disjoint
convex sets $C_1, C_2$ there is a halfspace $A$ so that
$C_1\subseteq A$ and $C_2\subseteq A^c$. See \cite[\S2]{Ro} for a
proof. The walls of a median algebra are pairs of complementary
halfspaces, so every two distinct points can be separated by
walls. However, the number of walls separating two distinct
points may be infinite, as it is the case with $\R$-trees.

We now show that discrete median algebras can be viewed as spaces
with walls and, correspondingly, morphisms of discrete median
algebras can be viewed as morphisms of spaces with walls
(Prop.3.7, equivalence of $a)$ and $c)$).

\begin{prop}
Any discrete median algebra is a space with walls.
\end{prop}
\proof
We need to show that, for every $x,y\in X$, the set $W(x,y)$ of
walls separating $x$ and $y$ is finite. We argue by induction that
$\big|W(x,y)\big|\leq \big|[x,y]\big|-1$.

Suppose $\big|[x,y]\big|=2$, i.e., $[x,y]=\{x,y\}$. Assume there
are two walls separating $x$ and $y$: $x\in A$, $y\in A^c$ and
$x\in B$, $y\in B^c$. Say $B\nsubseteq A$ and let $z\in
B\setminus A$. Then either $m(x,y,z)=x$, yielding $x\in A^c$, or
$m(x,y,z)=y$, yielding $y\in B$. Both are impossible.

Suppose $\big|[x,y]\big|>2$ and let $z\in [x,y]\setminus\{x,y\}$.
Then $[x,z]\cup [z,y]\subseteq [x,y]$ and $[x,z]\cap [z,y]$
consists of a single point, $z$, as $[x,z]\cap [z,y]=[x,z]\cap
[z,y]\cap [x,y]$. Using the induction step we obtain:
$$\big|W(x,y)\big|\leq \big|W(x,z)\big|+
\big|W(z,y)\big|\leq\big|[x,z]\big|-1+\big|[z,y]\big|-1\leq
\big|[x,y]\big|-1\eqno{\qed}$$

\begin{prop}
Let $f\co\negthinspace X\negthinspace\to\negthinspace X'$ be a
map, where $X$ and $X'$ are median algebras. The following are
equivalent:

{\rm(a)}\qua $f([x,y])\subseteq [f(x),f(y)]$ for all $x,y\in X$

{\rm(b)}\qua $f^{-1}(C')$ is convex in $X$ whenever $C'$ is convex in $X'$

{\rm(c)}\qua $f^{-1}(A')$ is a halfspace in $X$ whenever $A'$ is a
halfspace in $X'$

{\rm(d)}\qua $f\big(m(x,y,z)\big)=m\big(f(x),f(y),f(z)\big)$ for all
$x,y,z\in X$ \end{prop}

\begin{proof} This is a straightforward exercise in
median reasoning.

(a) $\Rightarrow$ (b):\qua Follows from the definition of convexity.

(b) $\Rightarrow$ (c):\qua Apply (b) to both $A'$ and $X'\setminus
A'$.

(c) $\Rightarrow$ (d):\qua Assume $f\big(m(x,y,z)\big)\neq
m\big(f(x),f(y),f(z)\big)$ for some $x,y,z\in X$. Then
$f\big(m(x,y,z)\big)\in A'$ and $m\big(f(x),f(y),f(z)\big)\in
X'\setminus A'$ for some halfspace $A'$ in $X'$. The latter
implies, by the convexity of $A'$, that at least two of
$\{f(x),f(y),f(z)\}$, say $f(x)$ and $f(y)$, are in $X'\setminus
A'$ i.e. $x,y\in f^{-1}(X'\setminus A')$. Then $m(x,y,z)\in
f^{-1}(X'\setminus A')$ as well, which is a contradiction.

(d) $\Rightarrow$ (a):\qua Let $z\in [x,y]$, that is $m(x,y,z)=z$.
Then $f(z)\in [f(x),f(y)]$ since
$m\big(f(x),f(y),f(z)\big)=f\big(m(x,y,z)\big)=f(z)$.
\end{proof}

The vertex set of a median graph, equipped with the interval
structure given by the geodesic intervals, is a discrete median
algebra. The converse holds as well: every discrete median algebra
arises as the vertex set of a median graph. Thus discrete median
algebras are precisely the 0--skeletons of median graphs. This
fact, which appears as a special case of our construction
(Cor.4.10), relates Chepoi's result \cite{Che00} that the
1--skeletons of $\cat$ cube complexes are precisely the median
graphs, to Roller's result \cite[Thm.10.3]{Ro}, that the
0--skeletons of $\cat$ cube complexes are precisely the discrete
median algebras.

\section{From spaces with walls to median graphs}
\noindent We are now ready to prove the main result:
\begin{thm}
Let $X$ be a space with walls. There exists an injective morphism
of spaces with walls
\begin{displaymath}
\begin{CD}
X @>\sigma>> \mathcal{C}^1(X)
\end{CD}
\end{displaymath}
where $\mathcal{C}^1(X)$ is a connected median graph, and
$\sigma(X)$ ``spans'' $\mathcal{C}^1(X)$, in the sense that no
proper subgraph of $\mathcal{C}^1(X)$ containing $\sigma(X)$ is
median.
\end{thm}

\noindent The core idea, which arises naturally in a variety of
contexts, can be summarized as follows: given a space, we identify
its points with the principal ultrafilters and then we suitably
add other ultrafilters. In our case, the suitable ultrafilters are
the almost principal ultrafilters. Here are the precise details.

\begin{defn}
Let $X$ be a space with walls. An \emph{ultrafilter} on $X$ is a
nonempty collection $\omega$ of halfspaces that satisfies the
following conditions:

$\bullet$\qua $A\in\omega$ and $A\subseteq B$ imply $B\in\omega$

$\bullet$\qua either $A\in\omega$ or $A^c\in\omega$ but not both
\end{defn}

Intuitively, an ultrafilter is a coherent orientation of the
walls. Note that every ultrafilter contains $X$. The easiest to
single out are the principal ultrafilters, defined for every $x\in
X$ to be the collection $\sigma_x$ of halfspaces containing $x$.

If $\omega_1$, $\omega_2$ are ultrafilters then elements of the
symmetric difference $\omega_1\triangle\omega_2$ come in pairs
$\{A,A^c\}$, so we may think of them as being walls. For distinct
$x$ and $y$, the set $\sigma_x\triangle\sigma_y$ describes the
walls separating $x$ and $y$, hence it is finite and nonempty.

Consider the graph $\Gamma$ whose vertices are the ultrafilters on
$X$, and the edges are defined by: $\omega_1$ is adjacent to
$\omega_2$ if $\frac{1}{2}\big|\omega_1\triangle\omega_2\big|=1$
i.e. $\omega_1$ and $\omega_2$ differ by exactly a wall. In
general $\Gamma$ is highly disconnected. The connectivity of
$\Gamma$ is described in the following statement.

\begin{lem}
There is a path in $\Gamma$ connecting $\omega_1$ and $\omega_2$
iff $\omega_1\triangle\omega_2$ is finite. In fact, the distance
between $\omega_1$ and $\omega_2$ is
$\frac{1}{2}\big|\omega_1\triangle\omega_2\big|$, which is the
number of walls by which $\omega_1$ and $\omega_2$ differ.
\end{lem}

\begin{proof}
Suppose $\omega_1=\theta_1,\dots,\theta_{m+1}=\omega_2$ is a path
connecting $\omega_1$ to $\omega_2$. Then
\begin{displaymath}
\frac{1}{2}\big|\omega_1\triangle\omega_2\big|=\frac{1}{2}\big|(\theta_1\triangle\theta_2)\triangle\dots\triangle(\theta_m\triangle\theta_{m+1})\big|\leq
\sum_{i=1}^m \frac{1}{2}\big|\theta_i\triangle\theta_{i+1}\big|=m
\end{displaymath}
Conversely, suppose $\omega_1\triangle\omega_2$ is finite. Let
$\omega_1\triangle\omega_2=\{A_1,\dots,A_n, A^c_1, \dots,
A^c_n\}$ where $A_i\in\omega_1\setminus\omega_2$ and
$A^c_i\in\omega_2\setminus\omega_1$. We may assume that each $A_i$
is minimal in $\{A_i,\dots,A_n\}$, and we define
$\theta_1=\omega_1$, $\theta_{i+1}=\theta_i\triangle\{A_i,
A^c_i\}$ for $1\leq i\leq n$. Note that $\theta_{n+1}=\omega_2$.

We claim that each $\theta_i$ is an ultrafilter. Since
$\theta_{i+1}$ is obtained from $\theta_i$ by exchanging $A_i$
for $A^c_i$, and since exchanging a minimal halfspace in an
ultrafilter for its complement results in an ultrafilter, we are
left with showing that $A_i$ is minimal in $\theta_i$. Suppose
there is $B\in\theta_i$, $B\subsetneq A_i$. Then
$B\notin\omega_2$ because $A_i\notin\omega_2$. As
\begin{displaymath}
\theta_i=\big(\omega_1\setminus\{A_1,\dots,A_{i-1}\}\big)\cup\{A^c_1,\dots,A^c_{i-1}\}
\end{displaymath}
we necessarily have $B\in\omega_1\setminus\{A_1,\dots,A_{i-1}\}$.
We obtain $B\in\{A_i,\dots,A_n\}$ which contradicts the fact that
$A_i$ is minimal in $\{A_i,\dots,A_n\}$.
\end{proof}

It follows that the principal ultrafilters lie in the same
connected component of $\Gamma$, denoted $\mathcal{C}^1(X)$. The
vertices of $\mathcal{C}^1(X)$ are the ultrafilters $\omega$ for
which $\omega\triangle\sigma_x$ is finite for some (every)
principal ultrafilter $\sigma_x$. We call them \emph{almost
principal ultrafilters}.

A helpful description of the geodesic intervals in
$\mathcal{C}^1(X)$ is the following:

\begin{lem}
Let $\omega, \omega_1,\omega_2$ be almost principal ultrafilters.
Then:
\begin{displaymath}
\omega\in[\omega_1,\omega_2]\Leftrightarrow\omega_1\cap\omega_2\subseteq\omega
\Leftrightarrow\omega\subseteq\omega_1\cup\omega_2
\end{displaymath}
\end{lem}

\begin{proof}
Since $\omega_1\triangle\omega$ and $\omega_2\triangle\omega$ are
finite, we have:
\begin{eqnarray*}
\omega\in[\omega_1,\omega_2]&\Leftrightarrow&
\big|\omega_1\triangle\omega\big|+\big|\omega\triangle\omega_2\big|=\big|\omega_1\triangle\omega_2\big|=\big|(\omega_1\triangle\omega)\triangle(\omega\triangle\omega_2)\big|\\
&\Leftrightarrow&
(\omega\triangle\omega_1)\cap(\omega\triangle\omega_2)=\emptyset\Leftrightarrow\omega_1\cap\omega_2\subseteq\omega\subseteq\omega_1\cup\omega_2
\end{eqnarray*}
The equivalence $\omega_1\cap\omega_2\subseteq\omega
\Leftrightarrow\omega\subseteq\omega_1\cup\omega_2$ holds for
arbitrary ultrafilters.
\end{proof}

\begin{prop}
$\mathcal{C}^1(X)$ is a median graph.
\end{prop}

\begin{proof}
Since geodesic intervals in $\mathcal{C}^1(X)$ are of the boolean
type described in Example~\ref{boolean}, the median in
$\mathcal{C}^1(X)$ of a triple of vertices $\omega_1$,
$\omega_2$, $\omega_3$, has to be the boolean median
\begin{displaymath}
m(\omega_1,\omega_2,\omega_3)=(\omega_1\cap\omega_2)\cup(\omega_2\cap\omega_3)\cup(\omega_3\cap\omega_1)\;
.
\end{displaymath}
We thus claim that $m(\omega_1,\omega_2,\omega_3)$ is a vertex in
$\mathcal{C}^1(X)$. Note that $m(\omega_1,\omega_2,\omega_3)$ is
an ultrafilter whenever $\omega_1$, $\omega_2$, $\omega_3$ are
ultrafilters. On the other hand, as $\omega_1\cap\omega_2\subseteq
m(\omega_1,\omega_2,\omega_3)\subseteq\omega_1\cup\omega_2$, we
have that $m(\omega_1,\omega_2,\omega_3)\triangle
\omega_2\subseteq\omega_1\triangle\omega_2$ so
$m(\omega_1,\omega_2,\omega_3)\triangle \omega_2$ is finite. Hence
$m(\omega_1,\omega_2,\omega_3)$ is almost principal.
\end{proof}

In order to show that the injective map $\sigma\co\negthinspace
X\negthinspace\to \mathcal{C}^1(X)$ given by $x\mapsto\sigma_x$
is the required embedding, we first need to understand the wall
structure of $\mathcal{C}^1(X)$.

\begin{prop}
There is a bijective correspondence between the halfspaces of $X$
and the halfspaces of $\mathcal{C}^1(X)$ given by
\begin{displaymath}
A\mapsto H_A=\{\omega\in\mathcal{C}^1(X):A\in\omega\}
\end{displaymath}
\end{prop}

\begin{proof}
Note that the complement of $H_A$ is $\mathcal{C}^1(X)\setminus
H_A=H_{A^c}$ and each $H_A$ is convex: if
$\omega\in[\omega_1,\omega_2]$ with $\omega_1,\omega_2\in H_A$,
then $A\in\omega_1\cap\omega_2\subseteq\omega$, hence $\omega\in
H_A$. Thus $H_A$ is a halfspace in $\mathcal{C}^1(X)$ for every
halfspace $A$ in $X$, which shows that the map is well-defined.

The map is injective since $\sigma_x$ is in $H_A$ iff $x\in A$,
i.e., $\sigma^{-1}(H_A)=A$.

We show that the map is surjective. Note that
$H_\emptyset=\emptyset$ and $H_X=\mathcal{C}^1(X)$. Let $H$ be a
proper halfspace in $\mathcal{C}^1(X)$. Pick $\omega\in H$,
$\omega'\notin H$ and consider a path
$\omega=\omega_0,\dots,\omega_n=\omega'$ connecting them. Then
$H$ cuts an edge in the path, in the sense that $\omega_i\in H$
and $\omega_{i+1}\notin H$ for some $i$. Suppose the edge
$\omega_i\omega_{i+1}$ is obtained by exchanging $A\in\omega_i$
for $A^c\in\omega_{i+1}$. We claim that $H=H_A$. If $\omega\in
H_A$, that is $A\in\omega$, then
$\omega_i\subseteq\omega_{i+1}\cup\omega$, i.e., $\omega_i\in
[\omega_{i+1},\omega]$. We get $\omega\in H$, since otherwise the
convexity of $\mathcal{C}^1(X)\setminus H$ would imply
$\omega_i\notin H$. Thus $H_A\subseteq H$. Similarly
$H_{A^c}\subseteq \mathcal{C}^1(X)\setminus H$, which by
complementation becomes $H\subseteq H_A$.
\end{proof}
We obtain a bijective correspondence between the walls of $X$ and
the walls of $\mathcal{C}^1(X)$ given by $\{A, A^c\}\mapsto
\{H_A,H_{A^c}\}$.

\begin{cor}
On $\mathcal{C}^1(X)$, the wall metric and the path metric
coincide.
\end{cor}
\begin{proof}
Let $\omega_1$, $\omega_2$ be vertices in $\mathcal{C}^1(X)$. The
wall metric counts the walls $\{H_A,H_{A^c}\}$ in
$\mathcal{C}^1(X)$ separating $\omega_1$, $\omega_2$. A wall
$\{H_A,H_{A^c}\}$ separates $\omega_1$, $\omega_2$ iff $\{A,
A^c\}\in\omega_1\triangle\omega_2$. The path metric counts the
walls $\{A, A^c\}$ in $\omega_1\triangle\omega_2$.
\end{proof}

\begin{prop}
The map $\sigma\co\negthinspace X\negthinspace \to
\negthinspace\mathcal{C}^1(X)$ given by $x\mapsto\sigma_x$ is an
injective morphism of spaces with walls, and an isometric
embedding when $X$ is equipped with the wall metric.
\end{prop}

\begin{proof} $\sigma$ is a morphism of spaces with walls
since each halfspace in $\mathcal{C}^1(X)$ is of the form $H_A$
for some halfspace $A$ in $X$, and $\sigma^{-1}(H_A)=A$. We have
that $d_w(x,y)=\frac{1}{2}\big|\sigma_x\triangle\sigma_y\big|$
and the right-hand side is the distance between $\sigma_x$ and
$\sigma_y$ in $\mathcal{C}^1(X)$.
\end{proof}

\begin{prop}
As a discrete median algebra, $\mathcal{C}^1(X)$ is generated by
the principal ultrafilters $\{\sigma_x:x\in X\}$.
\end{prop}
\begin{proof}
Let $M\subseteq\mathcal{C}^1(X)$ be the median subalgebra
generated by $\{\sigma_x:x\in X\}$. We proceed by contamination,
assuming that $\omega\in M$ and $\omega\omega '$ is an edge in
$\mathcal{C}^1(X)$, and proving that $\omega '\in M$. Suppose the
edge $\omega\omega '$ is obtained by exchanging $A^c\in\omega$
for $A\in\omega '$. Let $\zeta\in M\cap H_A$ be closest to
$\omega$. See Figure~\ref{fig:span}.

\begin{figure}[ht!] \centering
\includegraphics[width=.5\textwidth]{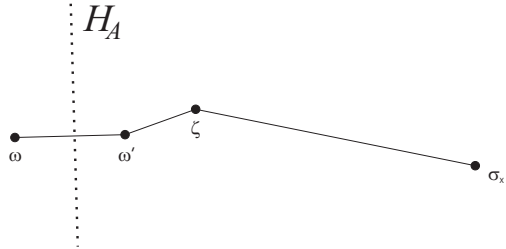}
\caption{The principal ultrafilters span} \label{fig:span}
\end{figure}

We claim that $\zeta=\omega '$. First, note that $\zeta\in
[\omega, \sigma_x]$ for all $x\in A$, since otherwise $m(\omega,
\zeta, \sigma_x)$ would be closer to $\omega$. Second, $\omega
'\in [\omega, \zeta]$. Hence $\zeta\in [\omega ', \sigma_x]$,
i.e., $\zeta\subseteq \omega ' \cup \sigma_x$ for all $x\in A$.
If there is $B\in\zeta\setminus\omega'$ then $B\in\sigma_x$ for
all $x\in A$, so $A\subseteq B$ and hence $B\in\omega '$, which
is a contradiction. Thus $\zeta\subseteq\omega'$, so
$\zeta=\omega'$.
\end{proof}

In particular, any discrete median algebra is isomorphic to a
median graph.
\begin{cor}
If $X$ is a discrete median algebra then $\sigma\co\negthinspace
X\negthinspace \to \negthinspace\mathcal{C}^1(X)$ is a median
isomorphism.
\end{cor}
\begin{proof}
One checks directly that
$m(\sigma_x,\sigma_y,\sigma_z)=\sigma_{m(x,y,z)}$. As
$\{\sigma_x:x\in X\}$ is closed under the median operation, it
equals $\mathcal{C}^1(X)$, in other words $\sigma$ is onto. Being
a bijective morphism of spaces with walls between median algebras,
$\sigma$ is a median isomorphism.
\end{proof}

Finally, we consider the problem of extending a group action on
$X$ to a group action on $\mathcal{C}^1(X)$.

\begin{prop}
Given a morphism of spaces with walls $f\co \negthinspace
X\negthinspace\to \negthinspace X'$, there exists a unique median
morphism $f_*\co \negthinspace\mathcal{C}^1(X)\negthinspace\to
\negthinspace\mathcal{C}^1(X')$ such that the following diagram
commutes:
\begin{displaymath}
\begin{CD}
X @>f>> X'\\
@V\sigma_X VV @VV\sigma_{X'} V\\
\mathcal{C}^1(X) @>{f_*}>> \mathcal{C}^1(X')
\end{CD}
\end{displaymath}
\end{prop}
\begin{proof}
Uniqueness is clear: a median morphism that makes the above
diagram commute is determined on $\{\sigma_x:x\in X\}$, which
spans $\mathcal{C}^1(X)$.

For the existence part, note that the condition
$f_*(\sigma_x)=\sigma_{f(x)}$ can be expressed as
$f_*(\sigma_x)=\{A'\subseteq X'\textrm{ halfspace
}:f^{-1}(A')\in\sigma_x\}$. This suggests the following
definition:
\begin{displaymath}
f_*(\omega)=\{A'\subseteq X'\textrm{ halfspace
}:f^{-1}(A')\in\omega\}
\end{displaymath}
Let us check that $f_*$ is well-defined. First, $f_*(\omega)$ is
an ultrafilter on $X'$ whenever $\omega$ is an ultrafilter on $X$.
Second, $f_*(\omega_1)\triangle f_*(\omega_2)$ is finite whenever
$\omega_1\triangle\omega_2$ is finite, since
\begin{displaymath}
f_*(\omega_1)\triangle f_*(\omega_2)=\{A'\subseteq X'\textrm{
halfspace }:f^{-1}(A')\in\omega_1\triangle\omega_2\}
\end{displaymath}
and the equation $f^{-1}(A')=A$, for a given proper halfspace $A$
of $X$, has finitely many solutions: picking $x\in A$, $y\in A^c$,
any solution separates $f(x)$ from $f(y)$.

\noindent To show that $f_*$ is a median morphism, one immediately
sees that $f_*$ preserves the boolean median. Alternatively, one
may show that $f_*$ is a morphism of spaces with walls. Indeed,
recall that every halfspace in $\mathcal{C}^1(X')$ is of the form
$H_{A'}=\{\omega '\in\mathcal{C}^1(X'):A'\in\omega '\}$ for some
halfspace $A'$ in $X'$. Then $f_*^{-1}(H_{A'})$ is a halfspace in
$\mathcal{C}^1(X)$:
\begin{eqnarray*}
f_*^{-1}(H_{A'})&=& \{\omega \in\mathcal{C}^1(X): f_*(\omega)\in
H_{A'}\}=\{\omega \in\mathcal{C}^1(X):A'\in
f_*(\omega)\}\\
&=&\{\omega\in\mathcal{C}^1(X):f^{-1}(A')\in\omega\}=H_{f^{-1}(A')}
\end{eqnarray*}

Note that $f_*$ need not be a graph morphism.\end{proof}

It follows that every automorphism of spaces with walls $f\co
\negthinspace X\negthinspace\to \negthinspace X$ has a unique
extension to a median automorphism, hence a graph automorphism as
well, $f_*\co \negthinspace\mathcal{C}^1(X)\negthinspace\to
\negthinspace\mathcal{C}^1(X)$ given by $f_*(\omega)=f(\omega)$.
Thus a group action on a space with walls $X$ naturally extends
to a group action on its 1--cubulation $\mathcal{C}^1(X)$.

\begin{figure}[ht!]\centering
\includegraphics[width=.5\textwidth]{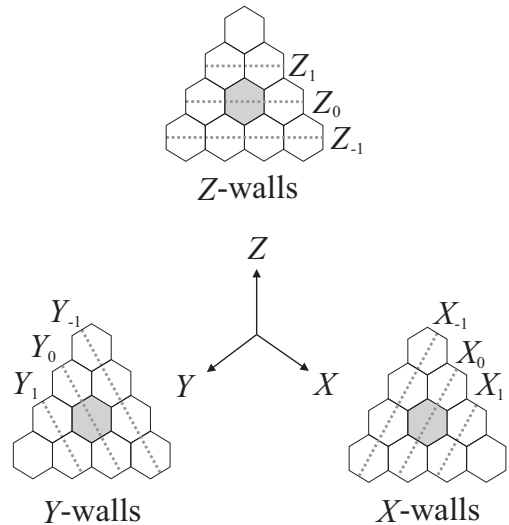}
\caption{Walls run across three directions \label{fig:hex}}
\end{figure}

\begin{ex} We cubulate the 1--skeleton of the hexagonal tiling of the plane.
The choice of halfspaces is independent along the three directions
$X$, $Y$, $Z$. But this is also the case with the choice of
halfspaces for the 1--skeleton of the usual tiling of $\R^3$ by
3--dimensional cubes. Since this is already a median graph, we
conclude that it is the 1--cubulation of the hexagonal tiling of
the plane.
\end{ex}

\Addresses\recd

\end{document}